\documentclass[11pt,bezier]{article}
\usepackage{amsmath}
\usepackage{amsfonts,amsthm,amssymb}
\usepackage{amsfonts}
\usepackage{graphicx}
\usepackage{epstopdf}
\usepackage{caption}
\usepackage{amssymb}
\newtheorem{lemma}{Lemma}[section]
\newtheorem{theorem}[lemma]{Theorem}

\usepackage{graphicx}
\theoremstyle{definition}

\theoremstyle{definition}

\usepackage[justification=centering]{caption}
\begin{document}
\begin{sloppypar}	

\title{The $g$-extra connectivity of the strong product of paths and cycles\footnote{The research is supported by National Natural Science Foundation of China (11861066).}}

\author{Qinze Zhu, Yingzhi Tian\footnote{Corresponding author. E-mail: qunzeZhu@qq.com (Q. Zhu); tianyzhxj@163.com (Y. Tian).} \\
{\small College of Mathematics and System Sciences, Xinjiang University, Urumqi, Xinjiang, 830046, PR China}}

\date{}
\maketitle
		
\noindent{\bf Abstract } Let $G$ be a connected graph and $g$ be a non-negative integer. The $g$-extra connectivity of $G$ is the minimum cardinality of a set of vertices in $G$, if it exists, whose removal disconnects $G$ and leaves every component with more than $g$ vertices. The strong product $G_1 \boxtimes G_2$ of graphs $G_1=(V_{1}, E_{1})$ and $G_2=(V_{2}, E_{2})$ is the graph with vertex set $V(G_1 \boxtimes G_2)=V_{1} \times V_{2}$, where two distinct vertices $(x_{1}, x_{2}), (y_{1}, y_{2}) \in V_{1} \times V_{2}$ are adjacent in $G_1 \boxtimes G_2$ if and only if $x_{i}=y_{i}$ or $x_{i} y_{i} \in E_{i}$ for $i=1, 2$. In this paper, we obtain the $g$-extra connectivity of the strong product of two paths, the strong product of a path and a cycle, and the strong product of two cycles.

	
\section{Introduction}
	
Let $G$ be a graph with vertex set $V(G)$ and edge set $E(G)$. The $minimum$ $degree$ of $G$ is denoted by $\delta(G)$. A $vertex$ $cut$ in $G$ is a set of vertices whose deletion makes $G$ disconnected. The $connectivity$ $\kappa(G)$ of the graph $G$ is the minimum order of a vertex cut in $G$ if $G$ is not  a complete graph; otherwise $\kappa(G)=|V(G)|-1$. Usually, the topology structure of an  interconnection network can be modeled by a graph $G$, where $V(G)$ represents the set of nodes and $E(G)$ represents the set of links connecting nodes in the network. Connectivity is used to measure the reliability the network, while it always underestimates the resilience of large networks.

To overcome this deficiency, Harary \cite{Harary} proposed the concept of conditional connectivity. For a graph-theoretic property $\mathcal{P}$, the $conditional$ $connectivity$ $\kappa(G; \mathcal{P})$ is the minimum cardinality of a set of vertices whose deletion disconnects $G$ and every remaining component has property $\mathcal{P}$. Later, F{\`{a}}brega and Fiol \cite{Fabrega} introduced the concept of $g$-extra connectivity, which is a kind of conditional connectivity. Let $g$ be a non-negative integer. A subset $S\subseteq V(G)$ is called a $g$-$extra$ $cut$ if $G-S$ is disconnected and each component of $G-S$  has at least $g+1$ vertices. The $g$-$extra$ $connectivity$ of $G$, denoted by $\kappa_{g}(G)$, is the minimum order of a $g$-extra cut if $G$ has at least one $g$-extra cut; otherwise define $\kappa_{g}(G)=\infty$. If $S$ is a $g$-extra cut in $G$ with order $\kappa_g(G)$, then we call $S$ a $\kappa_g$-$cut$. Since $\kappa_0(G)=\kappa(G)$ for any connected graph $G$ that is not a complete graph, the $g$-extra connectivity can be seen as a generalization of the traditional connectivity.  The authors in \cite{Chang} pointed out that there is  no polynomial-time algorithm for computing  $\kappa_g$ for a general graph. Consequently,  much of the work has been focused on the computing of the $g$-extra connectivity of some given graphs, see [1,4,6,8,10-11,16-21] for examples.

The most studied four standard graph products are the Cartesian product, the direct product, the strong product and the lexicographic product. The $Cartesian$ $product$ of two graphs $G_1$ and $G_2$, denoted by $G_1 \square G_2$, is defined on the vertex sets $V(G_1) \times V(G_2)$, and $(x_1, y_1)(x_2, y_2)$ is an edge in $G_1 \square G_2$ if and only if one of the following is true: ($i$) $x_1=x_2$ and $y_1 y_2 \in E(G_2)$;
($ii$) $y_1=y_2$ and $x_1 x_2 \in E(G_1)$.

The $strong$ $product$ $G_1 \boxtimes G_2$ of $G_1$ and $G_2$ is the graph with the vertex set $V(G_1 \boxtimes G_2)=V(G_1) \times V(G_2)$, where two vertices $(x_{1}, y_{1})$, $(x_{2}, y_{2}) \in V(G_1) \times V(G_2)$ are adjacent in $G_1 \boxtimes G_2$ if and only if one of the following holds: ($i$) $x_1=x_2$ and $y_1 y_2 \in E(G_2)$;
($ii$) $y_1=y_2$ and $x_1 x_2 \in E(G_1)$;
($iii$) $x_1 x_2 \in E(G_1)$ and $y_1 y_2 \in E(G_2)$.


{\v{S}}pacapan \cite{Spacapan1} proved that for any nontrivial graphs $G_1$ and $G_2$, $\kappa(G_1 \square G_2)=\min \{\kappa(G_1)|V(G_2)|, \kappa(G_2)|V(G_1)|, \delta(G_1 \square G_2)\}$.  L{\"{u}}, Wu,  Chen and Lv \cite{Lu} provided bounds for the 1-extra connectivity of the Cartesian product of two connected graphs.
Tian and Meng \cite{Tian} determined the exact values of the 1-extra connectivity of the Cartesian product for some class of graphs.  In \cite{Chen}, Chen, Meng, Tian and Liu further studied the 2-extra connectivity and the 3-extra connectivity of the Cartesian product of graphs.
	
Bre{\v{s}}ar and {\v{S}}pacapan \cite{Bresar} determined the edge-connectivity of the strong products of two connected graphs.
For the connectivity of the strong product graphs, {\v{S}}pacapan \cite{Spacapan2} obtained $Theorem\ \ref{2}$ in the following.
Let $S_i$ be a vertex cut in $G_i$ for $i=1,2$, and let $A_i$ be a component of $ G_i-S_i$ for $i=1,2$. Following the definitions in \cite{Spacapan2}, $I=S_1\times V_2$ or $I=V_1\times S_2$ is called an $I$-set in $G_1\boxtimes G_2$, and $L=(S_1\times A_2)\cup(S_1\times S_2)\cup(A_1\times S_2)$ is called an $L$-set in $G_1\boxtimes G_2$.
	
\begin{theorem}\label{2} (\cite{Spacapan2})
Let $G_1$ and $G_2$ be two connected graphs. Then every minimum vertex cut in $G_1\boxtimes G_2$ is either an $I$-set or an $L$-set in $G_1\boxtimes G_2$.
\end{theorem}
	
Motivated by the results above, we will study the $g$-extra connectivity of the strong product graphs. In the next section, we introduce some definitions and lemmas. In Section 3, we will give the $g$-extra connectivity of the strong product of two paths, the strong product of a path and a cycle, and the strong product of two cycles. Conclusion will be given in Section 4.

	\section{Preliminary}

For graph-theoretical terminology and notations not defined here, we follow \cite{Bondy}. Let $G$ be a graph with vertex set $V(G)$ and edge set $E(G)$. The $neighborhood$ of a vertex $u$ in $G$ is  $N_G(u)=\{v\in V(G)\ |\ v\;\text{is}\; \text{adjacent}\;\text{to} \;\text{the}\; \text{vertex}\; u\}$. Let $A$ be a subset of $V(G)$, the neighborhood of  $A$ in $G$ is $N_{G}(A)=\{v \in V(G) \backslash A \ |\ v\;\text{is}\; \text{adjacent}\;\text{to}\;  \text{a}\; \text{vertex}\; \text{in}\; A\}$. The subgraph induced by $A$ in $G$ is denoted by $G[A]$. We use $P_n$ to denote the path with order $n$ and $C_n$ to denote the cycle with order $n$.

Let  $G_1$ and $G_2$ be two graphs. Define two natural projections $p_1$ and $p_2$ on $V(G_1)\times V(G_2)$ as follows: $p_1(x,y)=x$ and $p_2(x,y)=y$ for any $(x,y)\in V(G_1)\times V(G_2)$.
The subgraph induced by $\{(u, y)|u\in V(G_1)\}$ in $G_1\boxtimes G_2$, denoted by $G_{1y}$, is called a $G_1$-$layer$ in $G_1\boxtimes G_2$ for each vertex $y\in V(G_2)$. Analogously, the subgraph induced by $\{(x, v)|v\in V(G_2)\}$  in $G_1\boxtimes G_2$, denoted by ${}_{x}G_2$, is called a $G_2$-$layer$ in $G_1\boxtimes G_2$ for each vertex $x\in V(G_1)$. Clearly,  a $G_1$-layer in $G_1\boxtimes G_2$ is isomorphic to $G_1$, and  a $G_2$-layer in $G_1\boxtimes G_2$ is isomorphic to $G_2$.
	
Let $S\subseteq V(G_1\boxtimes G_2)$. For any $x\in V(G_1)$, denote $S\cap V({}_{x}G_2)$ by ${}_{x}S$, and analogously,  for any $y\in V(G_2)$, denote $S\cap V(G_{1y})$ by $S_{y}$. Furthermore, we use $\overline{{}_{x}S}=V({}_{x}G_{2})\setminus {}_{x}S$ and $\overline{S_y}=V(G_{1y})\setminus S_y$. By a similar argument as the proof of the second paragraph of $Theorem$ 3.2 in \cite{Spacapan2}, we can obtain the following lemma.
	
	\begin{lemma}\label{1}
	Let $G$ be the strong product $G_1\boxtimes G_2$ of two connected graphs $G_1$ and $G_2$, and let $g$ be a non-negative integer. Assume $G$ has $g$-extra cuts and $S$ is a $\kappa_g$-cut of $G$.
		
(i) If ${}_{x}S\neq\emptyset$ for some $x\in V(G_1)$, then $|{}_{x}S|\geq \kappa(G_2)$.
		
(ii) If $S_{y}\neq \emptyset$ for some $y\in V(G_2)$, then $|S_{y}|\geq \kappa(G_1)$.
	
	\end{lemma}

	\noindent{\bf Proof.} ($i$) Suppose ${ }_{x} S \neq \emptyset$ for some $x\in V(G_1)$. Note that this is obviously true if ${ }_{x} S=V({ }_{x} G_{2})$. If $\overline{{}_{x}S}$ is not contained in one component of $G-S$, then clearly the induced subgraph $G[\overline{{}_{x}S}]$ is not connected, and hence $|{ }_{x} S| \geq \kappa(G_{2})$. If $\overline{{}_{x}S}$ is contained in one component of $G-S$, then choose an arbitrary fixed vertex $(x, y)$ from ${}_{x}S$. Let $H_{1}$ be the  component of $G-S$ such that $\overline{{}_{x}S} \subseteq V(H_{1})$ and let $H_{2}=G-S\cup V(H_{1})$. Since $S$ is a $\kappa_{g}$-cut, we find that the vertex $(x, y) \in{ }_{x} S$ has a neighbor $(x_{1}, y_{1}) \in V(H_{2})$. Since $(x_{1}, y_{1}) \in V(H_{2})$, we find that $(x, y_{1}) \in {}_{x} S$, moreover, for any $(x, u) \in \overline{{}_{x}S}$, we find that $(x, u)$ is not adjacent to $(x, y_{1})$, otherwise, $(x, u)$ would be adjacent to $(x_{1}, y_{1})$, which is not true since those two vertices are in different components of $G-S$. Thus if $R={ }_{x} S\setminus\{(x, y_{1})\}$, then $p_{2}(R)$ is a vertex cut in $G_{2}$ and one component of $G_{2}-p_{2}(R)$ is $\{y_{1}\}$. Thus $|{ }_{x} S|=|R|+1 \geq \kappa(G_{2})+1$.  Analogously, we can get $|S_{y}| \geq \kappa(G_{1})$ if ($ii$) holds.  $\Box$

	\section{Main results}	

	Let $H$ be a subgraph of $G_1\boxtimes G_2$. For the sake of simplicity, we use ${}_{x}H$ instead of ${}_{x}V(H)$ to represent $V(H)\cap V({}_{x}G_2)$ for any $x\in V(G_1)$ and $H_{y}$ to represent $V(H)\cap V(G_{1y})$ for any $y\in V(G_2)$.
Since  $\kappa_g(P_1\boxtimes P_n)=1$  for $g\leq \lfloor\frac{n-1}{2}\rfloor-1$ and $\kappa_g(P_2\boxtimes P_n)=2$  for $g\leq 2\lfloor\frac{n-1}{2}\rfloor-1$, we assume $m,n\geq3$ in the following theorem.

	\begin{theorem}\label{5}
		Let $g$ be a non-negative integer and $G=P_m\boxtimes P_n$, where $m,n\geq 3$. If $g\leq min\{n\lfloor\frac{m-1}{2}\rfloor-1, m\lfloor\frac{n-1}{2}\rfloor-1 \}$, then $\kappa_g(G)=min\{m, n, \lceil 2\sqrt{g+1}\ \rceil+1\}$.
	\end{theorem}

	\noindent{\bf Proof. }Denote $P_m=x_1x_2\dots x_m$ and $P_n=y_1y_2\dots y_n$.
Let $S_1=V(P_m)\times \{y_{\lfloor\frac{n-1}{2}\rfloor+1}\}$ and $S_2=\{x_{\lfloor\frac{m-1}{2}\rfloor+1}\}\times V(P_n)$.
Since $g\leq min\{n\lfloor\frac{m-1}{2}\rfloor-1, m\lfloor\frac{n-1}{2}\rfloor-1 \}$, we verify that $S_1$ and $S_2$ are two $g$-extra cuts of $G$. Thus $\kappa_g(G)\leq$ min$\{m, n\}$. If $\lceil 2\sqrt{g+1}\ \rceil+1\geq$ min$\{m, n\}$, then $\kappa_g(G)\leq$ min$\{m, n, \lceil 2\sqrt{g+1}\ \rceil+1\}$. If $\lceil 2\sqrt{g+1}\ \rceil+1<$ min$\{m, n\}$, then let $S_3=(J_1\times K_2)\cup (J_1\times J_2)\cup (K_1\times J_2)$, where
$J_1=\{x_{\lceil \sqrt{g+1}\ \rceil+1}\}$, $K_1=\{x_1, x_2, \dots, x_{\lceil \sqrt{g+1}\ \rceil}\}$, $J_2=\{y_{\lceil\frac{g+1}{\lceil \sqrt{g+1}\ \rceil}\rceil+1}\}$ and $K_2=\{y_1, y_2, \dots, y_{\lceil\frac{g+1}{\lceil \sqrt{g+1}\ \rceil}\rceil}\}$. It is routine to verify that $S_3$ is a $g$-extra cut of $G$. By $|S_3|=\lceil \sqrt{g+1}\ \rceil+\lceil\frac{g+1}{\lceil \sqrt{g+1}\ \rceil}\rceil+1=\lceil 2\sqrt{g+1}\ \rceil+1$,  we have $\kappa_g(G)\leq  \lceil 2\sqrt{g+1}\ \rceil+1$. Therefore, $\kappa_g(G)\leq$ min$\{m, n, \lceil 2\sqrt{g+1}\ \rceil+1\}$ holds.

Now, it is sufficient  to prove $\kappa_g(G)\geq$ min$\{m, n, \lceil 2\sqrt{g+1}\ \rceil+1\}$. Assume $S$ is a $\kappa_g$-cut of $G$. We consider two cases in the following.

	\noindent{\bf Case 1. } ${}_{x}S\neq\emptyset$ for all $x\in V(P_m)$, or $S_y\neq \emptyset$ for all $y\in V(P_n)$.
	
	Assume ${}_{x}S\neq\emptyset$ for all $x\in V(P_m)$. By $Lemma$ 2.1, $|S|= \sum_{x\in V(P_m)}|{}_{x}S|\geq \kappa(P_n)|V(P_m)|=m$. Analogously, if $S_y\neq \emptyset$ for all $y\in V(P_n)$, then $|S|=\sum_{y\in V(P_n)}|S_y|\geq \kappa(P_m)|V(P_n)|=n$.
	
	\noindent{\bf Case 2. } There exist a vertex $x_a\in V(P_m)$ and a vertex $y_b\in V(P_n)$   such that ${}_{x_a}S=S_{y_b}=\emptyset$.
	
By the assumption ${}_{x_a}S=S_{y_b}=\emptyset$, we know $V({}_{x_a}G_2)$  and $V(G_{1y_b})$ are contained in a component $H'$ of $G-S$. Let $H$ be another component of $G-S$. Let $p_1(V(H))=\{x_{s+1},x_{s+2},\cdots,x_{s+k}\}$ and $p_2(V(H))=\{y_{t+1},y_{t+2},\cdots,y_{t+h}\}$. Without loss of generality, assume $s+k<a$ and $t+h< b$. Clearly, $|V(H)|\leq kh$. Since $S$ is a $\kappa_g$-cut, we have $N_G(V(H))=S$ and $|V(H)|\geq g+1$. If we can prove $|N_G(V(H))|\geq k+h+1$, then $\kappa_g(G)=|S|=|N_G(V(H))|\geq k+h+1\geq2\sqrt{kh}+1\geq2\sqrt{g+1}+1$ and the theorem holds. Thus, we only need to show that $|N_G(V(H))|\geq k+h+1$ in the remaining proof.

\begin{figure}[htbp]
  \centering
  \includegraphics[width=12cm]{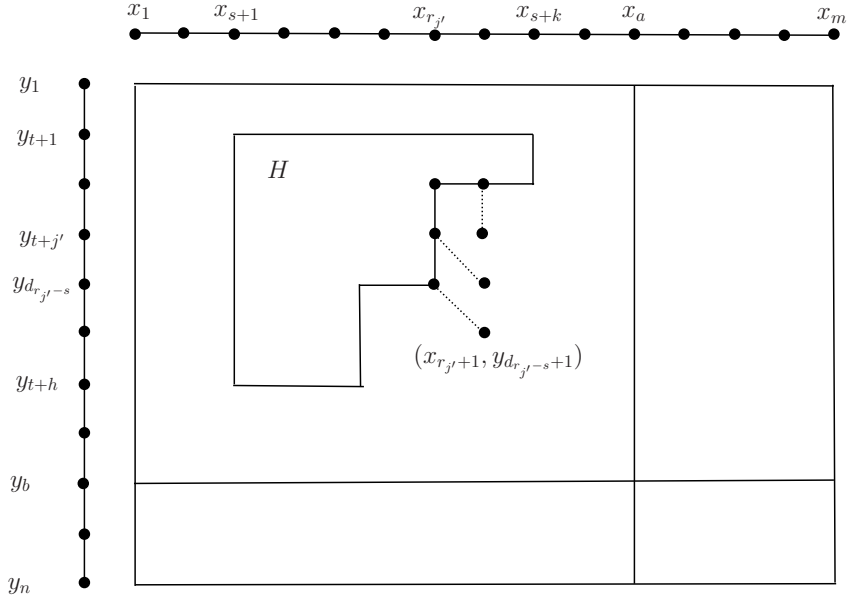}\\
  \caption{An illustration for the proof of Theorem 3.1.}
  \label{1}
\end{figure}

Let $(x_{s+i}, y_{d_i})$ be the vertex in ${}_{x_{s+i}}H$  such that $d_i$ is maximum for $i=1,\cdots, k$, and let $(x_{r_j}, y_{t+j})$ be the vertex in $H_{y_{t+j}}$  such that $r_j$ is maximum for $j=1,\cdots, h$.  Denote $D=\{(x_{s+1}, y_{d_1}),\cdots,(x_{s+k}, y_{d_k})\}$ and $R=\{(x_{r_1}, y_{t+1}),\cdots,(x_{r_h}, y_{t+h})\}$. For the convenience of counting, we will construct an injective mapping $f$ from $D\cup R$ to $N_G(V(H)\setminus\{(x_{s+k+1}, y_{d_k+1})\}$. Although $D$ and $R$ may have common elements, we consider the elements in $D$ and  $R$ to be different in defining the mapping $f$ below.

First, the mapping $f$ on $D$ is defined as follows.
\begin{center}
$f((x_{s+i}, y_{d_i}))=(x_{s+i}, y_{d_i+1})$ for $i=1,\cdots, k$.
\end{center}
Denote $F_1=\{(x_{s+1}, y_{d_1+1}),\cdots,(x_{s+k}, y_{d_k+1})\}$.

Second, for each vertex $(x_{r_j}, y_{t+j})$ satisfying $(x_{r_j+1}, y_{t+j})\notin F_1$, define $f((x_{r_j}, y_{t+j}))=(x_{r_j+1}, y_{t+j})$.

If $(x_{r_j}, y_{t+j})$ satisfies $(x_{r_j+1}, y_{t+j})\notin F_1$ for any $j\in\{1,\cdots, h\}$, then we are done.
Otherwise,  for each $(x_{r_{j'}}, y_{t+j'})$ satisfying $(x_{r_{j'}+1}, y_{t+j'})\in F_1$, we give the definition as follows.
By the definitions of $D$ and $R$, we have $(x_{r_{j'}+1+i}, y_{t+j'+j})\notin V(H)$ for all $i, j\geq 0$, and $\{(x_{r_{j'}},y_{t+j'}),\cdots,(x_{r_{j'}}, y_{d_{r_{j'}-s}})\}\subseteq R$ (see Figure 1 for an illustration). Now, we define $f((x_{r_{j'}}, y_{t+j'}))=(x_{r_{j'}+1}, y_{t+j'+1})$ and  change the images of $(x_{r_{j'}}, y_{t+j'+1}),\cdots,(x_{r_{j'}}, y_{d_{r_{j'}-s}})$ to $(x_{r_{j'}+1}, y_{t+j'+2}),\cdots,(x_{r_{j'}+1}, y_{d_{r_{j'}-s}+1})$, respectively.  The images of $f$ on $R$ are  well defined.

Finally, we have  an injective mapping $f$ from $D\cup R$ to $N_G(V(H)\setminus\{(x_{s+k+1}, y_{d_k+1})\}$. Then $\kappa_g(G)=|S|=|N_G(V(H))|\geq |D|+|R|+1\geq k+h+1\geq2\sqrt{kh}+1\geq2\sqrt{g+1}+1$. The proof is thus complete.
$\Box$

Since  $\kappa_g(C_3\boxtimes P_n)=3$  for $g\leq 3\lfloor\frac{n-1}{2}\rfloor-1$, we assume $m\geq4$ in the following theorem.

\begin{theorem}\label{6}
Let $g$ be a non-negative integer and $G=C_m\boxtimes P_n$, where $m\geq4$, $n\geq 3$. If $g\leq min\{n\lfloor\frac{m-2}{2}\rfloor-1, m\lfloor\frac{n-1}{2}\rfloor-1 \}$, then $\kappa_g(G)= min\{m, 2n, \lceil2\sqrt{2(g+1)}\rceil+2\}$.	
\end{theorem}

\noindent{\bf Proof.}
Denote $C_m=x_0x_1\dots x_{m-1} x_{m}$ (where $x_0=x_{m}$) and $P_n=y_1y_2\dots y_n$.
The addition of the subscripts of $x$ in the proof is modular $m$ arithmetic.
Let $S_1=V(C_m)\times\{y_{\lfloor\frac{n-1}{2}\rfloor+1}\}$ and $S_2=\{x_{0}, x_{\lfloor\frac{m-2}{2}\rfloor+1}\}\times V(P_n)$.
Since $g\leq min\{n\lfloor\frac{m-2}{2}\rfloor-1, m\lfloor\frac{n-1}{2}\rfloor-1 \}$, it is routine to check that $S_1$ and $S_2$ are two $g$-extra cuts of $G$. Thus $\kappa_g(G)\leq$ min$\{m, 2n\}$. If $\lceil2\sqrt{2(g+1)}\rceil+2\geq$ min$\{m, 2n\}$, then $\kappa_g(G)\leq$ min$\{m, 2n, \lceil2\sqrt{2(g+1)}\rceil+2\}$. If $\lceil2\sqrt{2(g+1)}\rceil+2<$ min$\{m, 2n\}$, then let $S_3=(J_1\times K_2)\cup (J_1\times J_2)\cup (K_1\times J_2)$, where
$J_1=\{x_{0}, x_{\lceil\sqrt{2(g+1)}\rceil+1} \}$, $K_1=\{x_1, x_2, \dots, x_{\lceil\sqrt{2(g+1)}\rceil}\}$, $J_2=\{y_{\lceil\frac{2(g+1)}{\lceil \sqrt{2(g+1)}\rceil}\rceil+1}\}$ and $K_2=\{y_1, y_2, \dots, y_{\lceil\frac{2(g+1)}{\lceil \sqrt{2(g+1)}\rceil}\rceil}\}$. It is routine to verify that $S_3$ is a $g$-extra cut of $G$. By $|S_3|=\lceil \sqrt{2(g+1)}\rceil+\lceil\frac{2(g+1)}{\lceil \sqrt{2(g+1)}\rceil}\rceil+2=\lceil 2\sqrt{2(g+1)}\rceil+2$,  we have $\kappa_g(G)\leq  \lceil2\sqrt{2(g+1)}\rceil+2$. Therefore, $\kappa_g(G)\leq$ min$\{m, 2n, \lceil2\sqrt{2(g+1)}\rceil+2\}$.
	
Now, it is sufficient  to prove $\kappa_g(G)\geq$ min$\{m, 2n, \lceil2\sqrt{2(g+1)}\rceil+2\}$. Assume $S$ is a $\kappa_g$-cut of $G$. We consider two cases in the following.

	\noindent{\bf Case 1. } ${}_{x}S\neq\emptyset$ for all $x\in V(C_m)$, or $S_y\neq \emptyset$ for all $y\in V(P_n)$.
	
	Assume ${}_{x}S\neq\emptyset$ for all $x\in V(C_m)$. By $Lemma$ 2.1, $|S|= \sum_{x\in V(C_m)}|{}_{x}S|\geq \kappa(P_n)|V(C_m)|=m$. Analogously, if $S_y\neq \emptyset$ for all $y\in V(P_n)$, then $|S|=\sum_{y\in V(P_n)}|S_y|\geq \kappa(C_m)|V(P_n)|=2n$.
	
	\noindent{\bf Case 2. } There exist a vertex $x_a\in V(C_m)$ and a vertex $y_b\in V(P_n)$   such that ${}_{x_a}S=S_{y_b}=\emptyset$.
	
By the assumption ${}_{x_a}S=S_{y_b}=\emptyset$, we know $V({}_{x_a}G_2)$  and $V(G_{1y_b})$ are contained in a component $H'$ of $G-S$. Let $H$ be another component of $G-S$. Let $p_1(V(H))=\{x_{s+1},x_{s+2},\cdots,x_{s+k}\}$ and $p_2(V(H))=\{y_{t+1},y_{t+2},\cdots,y_{t+h}\}$. Without loss of generality, assume $s+k<a$ and $t+h< b$. Clearly, $|V(H)|\leq kh$. Since $S$ is a $\kappa_g$-cut, we have $N_G(V(H))=S$ and $|V(H)|\geq g+1$. If we can prove $|N_G(V(H))|\geq k+2h+2$, then $\kappa_g(G)=|S|=|N_G(V(H))|\geq k+2h+2\geq2\sqrt{2kh}+2\geq2\sqrt{2(g+1)}+2$ and the theorem holds. Thus, we only need to show that $|N_G(V(H))|\geq k+2h+2$ in the remaining proof.

Let $(x_{s+i}, y_{d_i})$ be the vertex in ${}_{x_{s+i}}H$  such that $d_i$ is maximum for $i=1,\cdots, k$, and let $(x_{l_j}, y_{t+j})$ and $(x_{r_j}, y_{t+j})$ be the vertices in $H_{y_{t+j}}$  such that   $l_j$  and  $r_j$  are listed in the foremost and in the last along the sequence $(a+1,\cdots,m-1,0,1,\cdots,a-1)$, respectively, for $j=1,\cdots, h$.  Denote $D=\{(x_{s+1}, y_{d_1}),\cdots,(x_{s+k}, y_{d_k})\}$, $L=\{(x_{l_1}, y_{t+1}),\cdots,(x_{l_h}, y_{t+h})\}$ and $R=\{(x_{r_1}, y_{t+1}),\cdots,(x_{r_h}, y_{t+h})\}$. For the convenience of counting, we will construct an injective mapping $f$ from $D\cup L\cup R$ to $N_G(V(H)\setminus\{(x_{s}, y_{d_1+1}),(x_{s+k+1}, y_{d_k+1})\}$. Although $D$, $L$ and $R$ may have common elements, we consider the elements in $D$, $L$ and $R$ to be different in defining the mapping $f$ below.

First, the mapping $f$ on $D$ is defined as follows.
\begin{center}
$f((x_{s+i}, y_{d_i}))=(x_{s+i}, y_{d_i+1})$ for $i=1,\cdots, k$.
\end{center}
Denote $F_1=\{(x_{s+1}, y_{d_1+1}),\cdots,(x_{s+k}, y_{d_k+1})\}$.

Second, for each vertex $(x_{r_j}, y_{t+j})$ satisfying $(x_{r_j+1}, y_{t+j})\notin F_1$, define $f((x_{r_j}, y_{t+j}))=(x_{r_j+1}, y_{t+j})$.

If $(x_{r_j}, y_{t+j})$ satisfies $(x_{r_j+1}, y_{t+j})\notin F_1$ for any $j\in\{1,\cdots, h\}$, then we are done.
Otherwise,  for each $(x_{r_{j'}}, y_{t+j'})$ satisfying $(x_{r_{j'}+1}, y_{t+j'})\in F_1$, we give the definition as follows.
By the definitions of $D$ and $R$, we have  $\{(x_{r_{j'}},y_{t+j'}),\cdots,(x_{r_{j'}}, y_{d_{r_{j'}-s}})\}\subseteq R$. Now, we define $f((x_{r_{j'}}, y_{t+j'}))=(x_{r_{j'}+1}, y_{t+j'+1})$ and  change the images of $(x_{r_{j'}}, y_{t+j'+1}),\cdots,(x_{r_{j'}}, y_{d_{r_{j'}-s}})$ to $(x_{r_{j'}+1}, y_{t+j'+2}),\cdots,(x_{r_{j'}+1}, y_{d_{r_{j'}-s}+1})$, respectively.  The mapping $f$ on $R$ is defined well .

Third, for each vertex $(x_{l_j}, y_{t+j})$ satisfying $(x_{l_j-1}, y_{t+j})\notin F_1$, define $f((x_{l_j}, y_{t+j}))=(x_{l_j-1}, y_{t+j})$.

If $(x_{l_j}, y_{t+j})$ satisfies $(x_{l_j-1}, y_{t+j})\notin F_1$ for any $j\in\{1,\cdots, h\}$, then we are done.
Otherwise, for each $(x_{l_{j'}}, y_{t+j'})$ satisfying  $(x_{l_{j'}-1}, y_{t+j'})\in F_1$. By the definitions of $D$ and $L$, we have  $\{(x_{l_{j'}},y_{t+j'}),\cdots,(x_{l_{j'}}, y_{d_{l_{j'}-s}})\}\subseteq L$. Now, we define $f((x_{l_{j'}}, y_{t+j'}))=(x_{l_{j'}-1}, y_{t+j'+1})$ and  change the images of $(x_{l_{j'}}, y_{t+j'+1}),\cdots,(x_{l_{j'}}, y_{d_{l_{j'}-s}})$ to $(x_{r_{j'}-1}, y_{t+j'+2}),\cdots,(x_{l_{j'}-1}, y_{d_{r_{j'}-s}+1})$, respectively. The definition of $f$ on $L$ is complete.

Finally, we construct an injective mapping $f$ from $D\cup L\cup R$ to $N_G(V(H)\setminus\{(x_{s}, y_{d_1+1}),(x_{s+k+1}, y_{d_k+1})\}$. Then $\kappa_g(G)=|S|=|N_G(V(H))|\geq |D|+|L|+|R|+2\geq k+2h+1\geq2\sqrt{2kh}+2\geq2\sqrt{2(g+1)}+2$. The proof is thus complete.
$\Box$

Since  $\kappa_g(C_3\boxtimes C_n)=6$  for $g\leq 3\lfloor\frac{n-2}{2}\rfloor-1$, we assume $m, n\geq4$ in the following theorem.

\begin{theorem}
Let $g$ be a non-negative integer and $G=C_m\boxtimes C_n$, where $m, n\geq 4$. If $g\leq min\{n\lfloor\frac{m-2}{2}\rfloor-1, m\lfloor\frac{n-2}{2}\rfloor-1 \}$, then $\kappa_g(G)=min\{2m, 2n, \lceil 4\sqrt{g+1}\ \rceil+4\}$.
\end{theorem}

\noindent{\bf Proof.}
Denote $C_m=x_0x_1\dots x_{m-1} x_{m}$ (where $x_0=x_{m}$) and $C_n=y_0y_1\dots y_n$ (where $y_0=y_{n}$). The addition of the subscripts of $x$ in the proof is modular $m$ arithmetic, and the addition of the subscripts of $y$ in the proof is modular $n$ arithmetic.
Let $S_1=V(C_m)\times\{y_{0}, y_{\lfloor\frac{n-2}{2}\rfloor+1}\}$ and $S_2=\{x_{0}, x_{\lfloor\frac{m-2}{2}\rfloor+1}\}\times V(C_n)$.
Since $g\leq min\{n\lfloor\frac{m-2}{2}\rfloor-1, m\lfloor\frac{n-2}{2}\rfloor-1 \}$, we can check that $S_1$ and $S_2$ are two $g$-extra cuts of $G$. Thus $\kappa_g(G)\leq$ min$\{2m, 2n\}$. If $\lceil 4\sqrt{g+1}\ \rceil+4\geq$ min$\{2m, 2n\}$, then $\kappa_g(G)\leq$ min$\{2m, 2n, \lceil 4\sqrt{g+1}\ \rceil+4\}$. If $\lceil 4\sqrt{g+1}\ \rceil+4<$ min$\{2m, 2n\}$, then let $S_3=(J_1\times K_2)\cup (J_1\times J_2)\cup (K_1\times J_2)$, where
$J_1=\{x_{0}, x_{\lceil \sqrt{g+1}\ \rceil+1}\}$, $K_1=\{x_1, x_2, \dots, x_{\lceil \sqrt{g+1}\ \rceil}\}$, $J_2=\{y_0,y_{\lceil\frac{g+1}{\lceil \sqrt{g+1}\ \rceil} \rceil+1}\}$ and $K_2=\{y_1, y_2, \dots, y_{\lceil\frac{g+1}{\lceil \sqrt{g+1}\ \rceil}\rceil}\}$. It is routine to verify that $S_3$ is a $g$-extra cut of $G$. By $|S_3|=2\lceil \sqrt{g+1}\rceil+2\lceil\frac{g+1}{\lceil \sqrt{g+1}\ \rceil}\rceil+4=\lceil 4\sqrt{g+1}\rceil+4$,  we have $\kappa_g(G)\leq  \lceil 4\sqrt{g+1}\ \rceil+4$. Therefore, $\kappa_g(G)\leq$ min$\{2m, 2n, \lceil 4\sqrt{g+1}\ \rceil+4\}$.

Now, it is sufficient  to prove $\kappa_g(G)\geq$ min$\{2m, 2n, \lceil 4\sqrt{g+1}\ \rceil+4\}$. Assume $S$ is a $\kappa_g$-cut of $G$. We consider two cases in the following.

	\noindent{\bf Case 1. } ${}_{x}S\neq\emptyset$ for all $x\in V(C_m)$, or $S_y\neq \emptyset$ for all $y\in V(C_n)$.
	
	Assume ${}_{x}S\neq\emptyset$ for all $x\in V(C_m)$. By $Lemma$ 2.1, $|S|= \sum_{x\in V(C_m)}|{}_{x}S|\geq \kappa(C_n)|V(C_m)|=2m$. Analogously, if $S_y\neq \emptyset$ for all $y\in V(C_n)$, then $|S|=\sum_{y\in V(C_n)}|S_y|\geq \kappa(C_m)|V(C_n)|=2n$.
	
	\noindent{\bf Case 2. } There exist a vertex $x_a\in V(C_m)$ and a vertex $y_b\in V(C_n)$   such that ${}_{x_a}S=S_{y_b}=\emptyset$.
	
By the assumption ${}_{x_a}S=S_{y_b}=\emptyset$, we know $V({}_{x_a}G_2)$  and $V(G_{1y_b})$ are contained in a component $H'$ of $G-S$. Let $H$ be another component of $G-S$. Let $p_1(V(H))=\{x_{s+1},x_{s+2},\cdots,x_{s+k}\}$ and $p_2(V(H))=\{y_{t+1},y_{t+2},\cdots,y_{t+h}\}$. Without loss of generality, assume $s+k<a$ and $t+h< b$. Clearly, $|V(H)|\leq kh$. Since $S$ is a $\kappa_g$-cut, we have $N_G(V(H))=S$ and $|V(H)|\geq g+1$. If we can prove $|N_G(V(H))|\geq 2k+2h+4$, then $\kappa_g(G)=|S|=|N_G(V(H))|\geq 2k+2h+4\geq4\sqrt{kh}+4\geq4\sqrt{g+1}+4$ and the theorem holds. Thus, we only need to show that $|N_G(V(H))|\geq 2k+2h+4$ in the remaining proof.

Let $(x_{s+i}, y_{t_i})$ and $(x_{s+i}, y_{d_i})$ be the vertices in ${}_{x_{s+i}}H$  such that  $t_i$  and  $d_i$  are listed in the foremost and in the last along the sequence $(b+1,\cdots,n-1,0,1,\cdots,b-1)$, respectively, for $i=1,\cdots, k$, and let $(x_{l_j}, y_{t+j})$ and $(x_{r_j}, y_{t+j})$ be the vertices in $H_{y_{t+j}}$  such that   $l_j$  and  $r_j$  are listed in the foremost and in the last along the sequence $(a+1,\cdots,m-1,0,1,\cdots,a-1)$, respectively, for $j=1,\cdots, h$.  Denote $D=\{(x_{s+1}, y_{d_1}),\cdots,(x_{s+k}, y_{d_k})\}$, $T=\{(x_{s+1}, y_{t_1}),\cdots,(x_{s+k}, y_{t_k})\}$, $L=\{(x_{l_1}, y_{t+1}),\cdots,(x_{l_h}, y_{t+h})\}$ and $R=\{(x_{r_1}, y_{t+1}),\cdots,(x_{r_h}, y_{t+h})\}$. For the convenience of counting, we will construct an injective mapping $f$ from $D\cup T\cup L\cup R$ to $N_G(V(H)\setminus\{(x_{s}, y_{d_1+1}), (x_{s}, y_{t_1-1}), (x_{s+k+1}, y_{d_k+1}), (x_{s+k+1}, y_{t_k-1})\}$. Although $D$, $T$, $L$ and $R$ may have common elements, we consider the elements in $D$, $T$, $L$ and $R$ to be different in defining the mapping $f$ below.

First, the mapping $f$ on $D$ is defined as follows.
\begin{center}
$f((x_{s+i}, y_{d_i}))=(x_{s+i}, y_{d_i+1})$ for $i=1,\cdots, k$.
\end{center}
Denote $F_1=\{(x_{s+1}, y_{d_1+1}),\cdots,(x_{s+k}, y_{d_k+1})\}$.

Second, the mapping $f$ on $T$ is defined as follows.
\begin{center}
$f((x_{s+i}, y_{t_i}))=(x_{s+i}, y_{t_i-1})$ for $i=1,\cdots, k$.
\end{center}
Denote $F_2=\{(x_{s+1}, y_{t_1-1}),\cdots,(x_{s+k}, y_{t_k-1})\}$.

Third, for each vertex $(x_{r_j}, y_{t+j})$ satisfying $(x_{r_j+1}, y_{t+j})\notin F_1$, define $f((x_{r_j}, y_{t+j}))=(x_{r_j+1}, y_{t+j})$.

If $(x_{r_j}, y_{t+j})$ satisfies $(x_{r_j+1}, y_{t+j})\notin F_1$ for any $j\in\{1,\cdots, h\}$, then we are done.
Otherwise, for each $(x_{r_{j'}}, y_{t+j'})$ satisfying  $(x_{r_{j'}+1}, y_{t+j'})\in F_1$, we define as follows.
By the definitions of $D$ and $R$, we have  $\{(x_{r_{j'}},y_{t+j'}),\cdots,(x_{r_{j'}}, y_{d_{r_{j'}-s}})\}\subseteq R$. Now, we define $f((x_{r_{j'}}, y_{t+j'}))=(x_{r_{j'}+1}, y_{t+j'+1})$ and  change the images of $(x_{r_{j'}}, y_{t+j'+1}),\cdots,(x_{r_{j'}}, y_{d_{r_{j'}-s}})$ to $(x_{r_{j'}+1}, y_{t+j'+2}),\cdots,(x_{r_{j'}+1}, y_{d_{r_{j'}-s}+1})$, respectively.

Fourth, for each vertex $(x_{r_j}, y_{t+j})$ satisfying $(x_{r_j+1}, y_{t+j})\notin F_2$, define $f((x_{r_j}, y_{t+j}))=(x_{r_j+1}, y_{t+j})$.

If $(x_{r_j}, y_{t+j})$ satisfies $(x_{r_j+1}, y_{t+j})\notin F_2$ for any $j\in\{1,\cdots, h\}$, then we are done.
Otherwise, for each $(x_{r_{j'}}, y_{t+j'})$ satisfying $(x_{r_{j'}+1}, y_{t+j'})\in F_2$, we define as follows.
By the definitions of $T$ and $R$, we have  $\{(x_{r_{j'}},y_{t+j'}),\cdots,(x_{r_{j'}}, y_{d_{t_{j'}-s}})\}\subseteq R$. Now, we define $f((x_{r_{j'}}, y_{t+j'}))=(x_{r_{j'}+1}, y_{t+j'-1})$ and  change the images of $(x_{r_{j'}}, y_{t+j'-1}),\cdots,(x_{r_{j'}}, y_{d_{t_{j'}-s}})$ to $(x_{r_{j'}+1}, y_{t+j'-2}),\cdots,(x_{r_{j'}+1}, y_{d_{r_{j'}-s}-1})$, respectively.

Note that the proof of four paragraphs above gives the definition of the mapping $f$ on $R$. In the following proof, we will give the definition of the mapping $f$ on $L$.

Fifth, for each vertex $(x_{l_j}, y_{t+j})$ satisfying  $(x_{l_j-1}, y_{t+j})\notin F_1$, define $f((x_{l_j}, y_{t+j}))=(x_{l_j-1}, y_{t+j})$.

If $(x_{l_j}, y_{t+j})$ satisfies $(x_{l_j-1}, y_{t+j})\notin F_1$ for any $j\in\{1,\cdots, h\}$, then we are done.
Otherwise, for each $(x_{l_{j'}}, y_{t+j'})$ satisfying $(x_{l_{j'}-1}, y_{t+j'})\in F_1$, we define as follows.
By the definitions of $D$ and $L$, we have  $\{(x_{l_{j'}},y_{t+j'}),\cdots,(x_{l_{j'}}, y_{d_{l_{j'}-s}})\}\subseteq L$. Now, we define $f((x_{l_{j'}}, y_{t+j'}))=(x_{l_{j'}-1}, y_{t+j'+1})$ and  change the images of $(x_{l_{j'}}, y_{t+j'+1}),\cdots,(x_{l_{j'}}, y_{d_{l_{j'}-s}})$ to $(x_{r_{j'}-1}, y_{t+j'+2}),\cdots,(x_{l_{j'}-1}, y_{d_{r_{j'}-s}+1})$, respectively.

Sixth, for each vertex $(x_{l_j}, y_{t+j})$ satisfying $(x_{l_j-1}, y_{t+j})\notin F_2$, define $f((x_{l_j}, y_{t+j}))=(x_{l_j-1}, y_{t+j})$.

If $(x_{l_j}, y_{t+j})$ satisfies $(x_{l_j-1}, y_{t+j})\notin F_2$ for any $j\in\{1,\cdots, h\}$, then we are done.
Otherwise, for each $(x_{l_{j'}}, y_{t+j'})$ satisfying any $(x_{t_{j'}-1}, y_{t+j'})\in F_2$, we define as follows.
By the definitions of $L$ and $T$, we have  $\{(x_{l_{j'}},y_{t+j'}),\cdots,(x_{l_{j'}}, y_{d_{t_{j'}-s}})\}\subseteq L$. Now, we define $f((x_{l_{j'}}, y_{t+j'}))=(x_{l_{j'}-1}, y_{t+j'-1})$ and  change the images of $(x_{l_{j'}}, y_{t+j'-1}),\cdots,(x_{l_{j'}}, y_{d_{t_{j'}-s}})$ to $(x_{l_{j'}-1}, y_{t+j'-2}),\cdots,(x_{l_{j'}-1}, y_{d_{t_{j'}-s}-1})$, respectively.

Finally, we construct an injective mapping $f$ from $D\cup T\cup L\cup R$ to $N_G(V(H)\setminus\{(x_{s}, y_{d_1+1}), (x_{s}, y_{t_1-1}), (x_{s+k+1}, y_{d_k+1}), (x_{s+k+1}, y_{t_k-1})\}$. Then $\kappa_g(G)=|S|=|N_G(V(H))|\geq |D|+|T|+|L|+|R|+4\geq 2k+2h+4\geq4\sqrt{kh}+4\geq4\sqrt{g+1}+4$. The proof is thus complete.
$\Box$

	\section{Conclusion}	
	
Graph products are used to construct large graphs from small ones. Strong product is one of the most studied four graph products. As a generalization of traditional connectivity, $g$-extra connectivity can be seen as a refined parameter to measure the reliability of interconnection networks. There is no polynomial-time algorithm to compute the $g\ (\geq1)$-extra connectivity for a general graph. In this paper, we determined  the $g$-extra connectivity of the strong product of two paths, the strong product of a path and a cycle,  and the strong product of two cycles. In the future work, we would like to investigate the  $g$-extra connectivity of the strong product of two general graphs.

%

\end{sloppypar}

\begin{thebibliography}{5}



\bibitem{Boesch} F. T. Boesch, Synthesis of reliable networks-a survey, IEEE Trans. Reliab. 35(3) (1986) 240-246.

\bibitem{Bondy} J. A. Bondy and U. S. R. Murty, Graph Theory, Graduate Texts in Mathematics 244, Springer, Berlin, 2008.

\bibitem{Bresar} B. Bre{\v{s}}ar, S. {\v{S}}pacapan, Edge-connectivity of strong products of graphs, Discuss. Math. Graph Theory 27(2) (2007) 333-343.

\bibitem{Chang} N.-W. Chang, C.-Y. Tsai, S.-Y. Hsieh, On 3-extra connectivity and 3-extra edge connectivity of folded hypercubes, IEEE Trans. Comput. 63(6) (2014) 1593-1599.

\bibitem{Chen} L. H. Chen, J. X. Meng, Y. Z. Tian, F. X. Liu, Restricted connectivity of Cartesian product graphs, IAENG Int. J. Appl. Math. 46(1) (2016) 58-63.

\bibitem{Esfahanian} A. Esfahanian, S. Hakimi, On computing a conditional edge-connectivity of a graph, Inform. Process. Lett. 27(4) (1988) 195-199.

\bibitem{Fabrega} J. F{\`{a}}brega, M. A. Fiol, On the extraconnectivity of graphs, Discrete Math. 155(1-3) (1996) 49-57.

\bibitem{Guo} H. M. Guo, E. Sabir, A. Mamut, The $g$-extra connectivity of folded crossed cubes, J. Parallel Distributed Comput. 166 (2022) 139-146.

\bibitem{Harary} F. Harary, Conditional connectivity, Networks 13(3) (1983) 347-357.

\bibitem{Hsieh} S.-Y. Hsieh, Y.-H. Chang, Extraconnectivity of $k$-ary $n$-cube networks, Theoret. Comput. Sci. 443(20) (2012) 63-69.

\bibitem{Latifi} S. Latifi, M. Hegde, M. Naraghi-Pour, Conditional connectivity measures for large multiprocessor systems, IEEE Trans. Comput. 43(2) (2002) 218-222.

\bibitem{Lu} M. L\"{u}, C. Wu, G.-L. Chen, C. Lv,  On super connectivity of Cartesian product graphs, Networks 52(2) (2008) 78-87.

\bibitem{Spacapan1} S. {\v{S}}pacapan, Connectivity of Cartesian products of graphs. Appl. Math. Lett. 21(7) (2008) 682-685.

\bibitem{Spacapan2}  S. {\v{S}}pacapan, Connectivity of Strong Products of Graphs, Graphs Comb. 26(3)(2010) 457-467.

\bibitem{Tian} Y. Z. Tian, J. X. Meng, Restricted connectivity for some interconnection networks, Graphs Comb. 31(5) (2015) 1727-1737.

\bibitem{Wan} M. Wan, Z. Zhang, A kind of conditional vertex connectivity of star graphs, Appl. Math. Lett. 22(2) (2009) 264-267.

\bibitem{Yang-Lin} W. H.  Yang, H. Q. Lin, Reliability evaluation of BC networks in terms of the extra vertex-and edge-connectivity, IEEE Trans. Comput. 63(10) (2014) 2540-2548.

\bibitem{Yang-Meng} W. H. Yang, J. X. Meng, Extraconnectivity of hypercubes, Appl. Math. Lett. 22(6) (2009) 887-891.

\bibitem{Zhang} M. M. Zhang, J. X. Zhou, On $g$-extra connectivity of folded hypercubes, Theor. Comput. Sci. 593 (2015) 146-153.

\bibitem{Zhou} J. X. Zhou, On $g$-extra connectivity of hypercube-like networks, J. Comput. Syst. Sci. 88 (2017) 208-219.

\bibitem{Zhu} Q. Zhu, X. K. Wang, G. L. Cheng, Reliability Evaluation of BC Networks, IEEE Trans. Computers 62(11) (2013) 2337-2340.
		
		
		
\end{thebibliography}
\end{document}